%% file: agt-2-17.tex
\documentclass{gtart}

\input agtout

\lognumber{17}
\volumenumber{2}
\volumeyear{2002}
\papernumber{17}
\published{21 May 2002}
\pagenumbers{371}{380}
\received{13 March 2002}
\accepted{28 March 2002}

\usepackage{amsmath, amssymb}
\usepackage{graphicx}

\input labelfig

\newtheorem{lemma}{Lemma}
\newtheorem*{thm1}{Theorem~1}
\newtheorem*{thm2}{Theorem~2}
\newtheorem*{thm3}{Theorem~3}

\begin{document}

\title {Intrinsic knotting and linking of complete graphs}

\author {Erica Flapan}

\address {Department of 
Mathematics, Pomona College\\
Claremont, CA 91711, U.S.A.}
\email{eflapan@pomona.edu}

\begin{abstract}We show that for 
every $m\in {\mathbb N}$, there exists an $n\in {\mathbb N}$ such that
every embedding of the complete graph $K_{n}$ in ${\mathbb R}^{3}$
contains a link of two components whose linking number is at least
$m$.  Furthermore, there exists an $r\in {\mathbb N}$ such that every
embedding of $K_{r}$ in ${\mathbb R}^{3}$ contains a knot $Q$ with $
\vert a_{2}(Q) \vert \geq m$, where $a_{2}(Q)$ denotes the second
coefficient of the Conway polynomial of $Q$.
\end{abstract}

\asciiabstract{We show that for every m in N, there exists an n in N
such that every embedding of the complete graph K_n in R^3 contains a
link of two components whose linking number is at least m.
Furthermore, there exists an r in N such that every embedding of K_r
in R^3 contains a knot Q with |a_2(Q)| > m-1, where a_2(Q) denotes the
second coefficient of the Conway polynomial of Q.}

\maketitle

\section{Introduction}   

The study of intrinsic knotting and linking began 
with the work of Conway and Gordon \cite{CG}, who 
showed that every embedding in ${\mathbb R}^{3}$ of 
the complete graph on six vertices, $K_{6}$, 
contains a non-trivial link of two components, 
and every embedding of $K_{7}$ in ${\mathbb R}^{3}$ 
contains a non-trivial knot.  Since the existence 
of such a link or knot is intrinsic to the graph, 
and does not depend on the particular embedding 
of the graph in ${\mathbb R}^{3}$, we say that 
$K_{6}$ is {\it intrinsically linked\/} and 
$K_{7}$ is {\it intrinsically knotted\/}.  We 
state Conway and Gordon's Theorem more precisely 
with the following notation.  Let $L_{1}\cup 
L_{2}$ be an oriented link, and let 
$lk(L_{1},L_{2})$ denote the linking number of 
$L_{1}$ and $L_{2}$; let $Q$ be a knot, and let 
$a_{2}(Q)$ denote the second coefficient of the 
Conway polynomial of $Q$.  Conway and Gordon 
proved that every embedding of $K_{6}$ contains a 
link $L_{1}\cup L_{2}$ which has 
$lk(L_{1},L_{2})\equiv 1$ (mod $2$), and every 
embedding of $K_{7}$ contains a knot $Q$ which 
has $a_{2}(Q)\equiv 1$ (mod $2$).  Furthermore, 
they illustrated an embedding of $K_{6}$ such 
that the only non-trivial link $L_{1}\cup L_{2}$ 
contained in $K_{6}$ is the Hopf link (which has 
$ \vert lk(L_{1},L_{2}) \vert =1$); and they 
illustrated an embedding of $K_{7}$ such that the 
only non-trivial knot $Q$ contained in $K_{7}$ is 
the trefoil knot (which has $ \vert a_{2}(Q) 
\vert =1$).  In this sense $K_{6}$ exhibits the 
simplest type of intrinsic linking and $K_{7}$ 
exhibits the simplest type of intrinsic knotting.

	In this paper, we will show that for larger 
values of $n$, the complete graph $K_{n}$ can 
exhibit a more complex type of intrinsic linking 
or knotting.  In particular, there exists an $n$ 
such that every embedding of $K_{n}$ contains a 
non-trivial $2$-component link which is not the 
Hopf link, and there exists an $n$ such that 
every embedding of $K_{n}$ contains a non-trivial 
knot which is not the trefoil knot.  In 
\cite{FNP} we considered links of more than two 
components, and showed that every embedding of 
$K_{10}$ in ${\mathbb R}^{3}$ contains a 
$3$-component link $L_{0}\cup L_{1}\cup L_{2}$ 
such that both $lk(L_{0},L_{1})\equiv 1$ (mod 
$2$) and $lk(L_{0},L_{2})\equiv 1$ (mod $2$); 
furthermore, $n=10$ is the smallest number such 
that $K_{n}$ has this property.  Here we shall 
generalize Conway and Gordon's Theorem by 
considering the linking number of 2-component 
links in ${\mathbb Z}$ rather than in ${\mathbb 
Z}_{2}$, and by considering the second 
coefficient of the Conway polynomial in ${\mathbb 
Z}$ rather than in ${\mathbb Z}_{2}$. 

We begin by proving that every embedding of 
$K_{10}$ contains a non-trivial link of two 
components other than the Hopf link.  In 
particular we prove the following.

\begin{thm1}
Every embedding
of 
$K_{10}$ in ${\mathbb R}^{3}$ contains a 
$2$-component link $L=L_{1}\cup J_{1}$ such that 
for some orientation of $L$ we have 
$lk(L_{1},J_{1})\geq 2$.  
\end{thm1}

Theorem~2 will show that the complexity of
the  intrinsic linking of $K_{n}$ (as measured by
the  linking number) can be made as large as we wish 
by making $n$ sufficiently large. 

\begin{thm2}
  Let $p\in {\mathbb N}$ be 
given, and let $n=p(15p-9)$.  Then every 
embedding of $K_{n}$ in ${\mathbb R}^{3}$ contains a 
$2$-component link $L=L_{p}\cup J_{p}$ such that 
for some orientation of $L$ we have 
$lk(L_{p},J_{p})\geq p$. 
\end{thm2}     

\begin{figure}[ht!]
\label{fig1}
\centerline{\small
%\ShowGrid 
\SetLabels 
\E(0.13*-.1){$K_+$}\\
\E(0.55*-.1){$K_-$}\\
\E(0.83*.5){$L_1$}\\
\E(0.99*.5){$L_2$}\\
\endSetLabels 
\AffixLabels{{\includegraphics[width=2.5in]{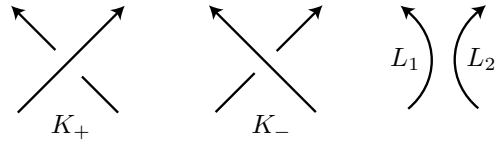}}}}
\vspace{2mm}
\caption{The skein moves}
\end{figure}

It is natural to ask whether the complexity of 
the intrinsic knotting of $K_{n}$ can also be 
made as large as we wish by choosing $n$ 
sufficiently large.  While the linking number 
seems like the natural measure of the complexity 
of a 2-component link, there are various ways to 
measure the complexity of a knot.  The second 
coefficient of the Conway polynomial $a_{2}(Q)$ 
of an oriented knot $Q$ is a convenient invariant 
to use because it relates knotting and linking.  
In particular, Kauffman \cite{Ka} has shown that 
it satisfies the following equation
\begin{equation}
\label{eqn1}
a_{2}(K_{+})=a_{2}(K_{-})+lk(L_{1},L_{2})   
\end{equation}
\par \noindent where $K_{+}$ and $K_{-}$ are 
identical oriented knots outside of the crossing 
illustrated in Figure~\ref{fig1}, and the oriented
link 
$L_{1}\cup L_{2}$ is obtained by smoothing this 
crossing as illustrated.

Due to the utility of equation $(1)$, the second 
coefficient of the Conway polynomial has been 
used to prove various theorems about knots 
contained in spatial graphs (see for example 
\cite{CG}, \cite{Fo}, \cite{TY}, \cite{Sh}).  
Using $a_{2}(Q)$ as a measure of knot complexity 
we prove the following theorem.

\begin{thm3}
 Let $m\in {\mathbb N}$ be 
given, let $p$ be an integer such that 
$p\geq 4\sqrt {m}$, and let $n=p(15p-9)$.  
Then every embedding of $K_{2n}$ in ${\mathbb 
R}^{3}$ contains a knot $Q$ with $ \vert a_{2}(Q) 
\vert \geq m$. 
\end{thm3}        

Observe that by picking $p$ to be the smallest 
integer such that $p\geq 4\sqrt {m}$, we 
have $p\leq 4\sqrt {m}+1$ and hence 
$2n\leq 720m-60$.  In particular, this means 
that $2n$ only grows linearly with $m$.

The question of finding the minimum number of 
vertices necessary to guarantee a certain type of 
intrinsic linking or knotting remains open.  
Since $K_{5}$ does not have enough vertices to 
contain a link, $n=6$ is the smallest number of 
vertices necessary such that every embedding of 
$K_{n}$ contains a link with $L_{1}\cup L_{2}$ 
with $ \vert lk(L_{1},L_{2}) \vert \geq 1$.  
Theorem 1 shows every embedding of $K_{10}$ in 
${\mathbb R}^{3}$ contains a link with $L_{1}\cup 
L_{2}$ with $ \vert lk(L_{1},L_{2}) \vert 
\geq 2$.   However, it is not known whether 
$10$ is the smallest such $n$.  In general, it is 
an open question to find a function $\varphi 
:{\mathbb N}\rightarrow {\mathbb N}$ such that for 
every $p\in {\mathbb N}$, $\varphi (p)$ is the 
minimum number of vertices necessary such that 
every embedding of $K_{\varphi (p)}$ contains a 
$2$-component link $L_{1}\cup L_{2}$ with $ \vert 
lk(L_{1},L_{2}) \vert \geq p$.  With respect 
to intrinsic knotting, the embedding of $K_{6}$ 
given by Conway and Gordon contains no 
non-trivial knot, thus $n=7$ is the minimum 
number of vertices necessary such that every 
embedding of $K_{n}$ contains a knot $Q$ with $ 
\vert a_{2}(Q) \vert \geq 1$.  It follows 
from Theorem 3 that every embedding of $K_{972}$ 
in ${\mathbb R}^{3}$ contains a non-trivial knot 
other than the trefoil knot.  However it is not 
known what the smallest $n$ is such that every 
embedding of $K_{n}$ in ${\mathbb R}^{3}$ contains a 
knot other than the trefoil knot.  Furthermore, 
it is an open question to find a function $\psi 
:{\mathbb N}\rightarrow {\mathbb N}$ such that for 
every $m\in {\mathbb N}$, $\psi (m)$ is the minimum 
number of vertices necessary such that every 
embedding of $K_{\psi (m)}$ contains a knot $Q$ 
with $ \vert a_{2}(Q) \vert \geq m$. 

\section{Intrinsic
Linking}

The following lemma allows us to go from a 
3-component link (with a sufficient number of 
vertices) to a 2-component link whose linking 
number is at least the sum of the linking numbers 
of two pairs of components of the 3-component 
link.
  
\begin{lemma}
\label{lem1}  Let $L\cup Z\cup W$ be a 
$3$-component link contained in some embedding of 
$K_{n}$ in ${\mathbb R}^{3}$.  Suppose that 
$lk(L,Z)=p_{1}>0$ and $lk(L,W)=p_{2}>0$ for some 
orientation of $L\cup Z\cup W$.  Suppose that $Z$ 
and $W$ each contain at least $q$ vertices and 
$q>p_{1}+p_{2}$.  Then $K_{n}$ contains a simple 
closed curve $J$ with at least $2q$ vertices which is
disjoint from $L$ such that, for some orientation of
$L\cup J$, we have 
$lk(L,J)\geq p_{1}+p_{2}$. 
\end{lemma}   
\begin{proof}   On $Z$ we select $q$ vertices and 
label them consecutively by $v_{1}$, \dots, $v_{q}$ 
in such a way that $Z$ is oriented in the 
direction of increasing order of the $v_{j}$.  On 
$W$ we select $q$ vertices and label them 
consecutively by $w_{1}$, \dots, $w_{q}$ in such a 
way that $W$ is oriented in the direction of 
decreasing order of the $w_{j}$.  For each 
$j=1$, \dots, $q$, let $A_{j}$ denote the simple 
closed curve $\overline 
{v_{j}w_{j}w_{j+1}v_{j+1}v_{j}}$, where the 
subscripts are taken mod $q$.  We orient each 
$A_{j}$ so that going from $v_{j}$ to $w_{j}$ 
along the edge $\overline {v_{j}w_{j}}$ is the 
positive direction.

Now in the homology group $H_{1}({\mathbb 
R}^{3}-L;{\mathbb Z})$ we have the equation
\begin{equation*}
[Z]+[W]+[A_{1}]+\ldots +[A_{q}]=0.
\end{equation*}
 Thus
$$
p_{1}+p_{2}=[Z]+[W]=-[A_{1}]- \ldots -[A_{q}].
$$
Since $[A_{j}]$ is an integer for each $j$ and 
$q>p_{1}+p_{2}$, there is some $j$ such that 
$[A_{j}]\geq 0$.  Without loss of generality 
$[A_{q}]\geq 0$.  Hence in $H_{1}({\mathbb 
R}^{3}-L;{\mathbb Z})$ we have the inequality 
$$
-[A_{1}]-\ldots -[A_{q-1}]\geq
p_{1}+p_{2}.
$$
Now let $J$ denote the simple closed curve 
obtained from $A_{1}\cup \ldots\cup A_{q-1}$ by 
omitting the edges $\overline {v_{j}w_{j}}$ for 
$j=2$, \dots, $q-1$.  We orient $J$ so that going 
from $w_{1}$ to $v_{1}$ along the edge $\overline 
{w_{1}v_{1}}$ is the positive direction.  Then in 
$H_{1}({\mathbb R}^{3}-L;{\mathbb Z})$ we have 
$$
[J]=-[A_{1}]-\ldots-[A_{q-1}].  
$$
Hence $lk(L,J)\geq p_{1}+p_{2}$ and $J$ has 
at least $2q$ vertices.  
\end{proof}

We shall prove Theorem~1 by using
Lemma~\ref{lem1}  together with \cite{FNP}.

\begin{thm1}   Every embedding of 
$K_{10}$ in ${\mathbb R}^{3}$ contains a 
$2$-component link $L\cup J$ such that, for some 
orientation, we have $lk(L,J)\geq 
2$.
\end{thm1}    
\begin{proof}   Let $K_{10}$ be embedded in 
${\mathbb R}^{3}$.  It follows from \cite{FNP} that 
$K_{10}$ contains a $3$-component link $L\cup 
Z\cup W$ such that $lk(L,Z)\equiv 1$ (mod $2$) 
and $lk(L,W)\equiv 1$ (mod $2$). We orient the 
link $L\cup Z\cup W$ so that 
$lk(L,Z)=p_{1}\geq 1$ and 
$lk(L,W)=p_{2}\geq 1$.  If either 
$p_{i}\geq 3$ we are done.  Otherwise, 
$p_{1}=p_{2}=1$.   Clearly $Z$ and $W$ each have 
at least $3$ vertices.  Now we apply Lemma 1 to 
obtain a simple closed curve $J$ such that 
$lk(L,J)\geq 2$.
\end{proof}

Lemma~\ref{lem2} will allow us to go from a pair of 
disjoint 2-component links (each with a 
sufficient number of vertices) to a 3-component 
link, while bounding the linking numbers of the 
3-component link below.  The proof of
Lemma~\ref{lem2} is  similar in flavor to the proof
of Lemma~\ref{lem1}, though  the details differ.

\begin{lemma}
\label{lem2}   Let $X_{1}\cup Y_{1}\cup 
X_{2}\cup Y_{2}$ be a $4$-component link 
contained in some embedding of $K_{n}$ in ${\mathbb 
R}^{3}$.  Suppose that for some orientation of 
$X_{1}\cup Y_{1}\cup X_{2}\cup Y_{2}$ we have 
$lk(X_{1},Y_{1})\geq 1$ and 
$lk(X_{2},Y_{2})=p\geq 1$.  Also suppose 
that $X_{1}$, $Y_{1}$, $X_{2}$, $Y_{2}$ each 
contain at least $q$ vertices and $q>p$.  Then 
$K_{n}$ contains disjoint simple closed curves $L$,
$Z$  and $W$, each with at least $q$ vertices, such 
that $lk(L,Z)\geq 1$ and $lk(L,W)\geq 
p$ for some orientation of $L\cup Z\cup W$.   
\end{lemma}  

\begin{proof}  If $lk(X_{2},Y_{1})$ is non-zero, 
then let $L=X_{2}$, $Z=Y_{1}$, and $W=Y_{2}$.  
Now if we orient $L\cup Z\cup W$ appropriately we 
get a link with $lk(L,Z)\geq 1$ and 
$lk(L,W)\geq p$.  If $lk(Y_{2},X_{1})$ is 
non-zero, let $L=Y_{2}$, $Z=X_{1}$, and 
$W=X_{2}$.  Then $L\cup Z\cup W$ is the desired 
link.  So from now on, we shall assume that both 
$lk(Y_{1},X_{2})=0$ and $lk(X_{1},Y_{2})=0$.

We choose $q$ vertices on each of $X_{1}$ and 
$X_{2}$ and label the vertices of $X_{1}$ and 
$X_{2}$ as we did the vertices of $Z$ and $W$ in 
the proof of Lemma 1. We also define the oriented 
simple closed curves $A_{j}$ as we did in the 
proof of Lemma 1.  Now for both $i=1$ and $i=2$, 
in the first homology groups $H_{1}({\mathbb 
R}^{3}-Y_{i};{\mathbb Z})$ we have the equation
$$
[X_{1}]+[X_{2}]+[A_{1}]+\ldots+[A_{q}]=0.
$$
Now by our assumptions that $lk(Y_{1},X_{2})=0$ 
and $lk(X_{1},Y_{2})=0$, we have $[X_{2}]=0$ in 
$H_{1}({\mathbb R}^{3}-Y_{1};{\mathbb Z})$ and we have 
$[X_{1}]=0$ in $H_{1}({\mathbb R}^{3}-Y_{2};{\mathbb 
Z})$.  Thus in $H_{1}({\mathbb R}^{3}-Y_{2};{\mathbb 
Z})$ we have
\begin{equation*}
0<p=[X_{2}]=-[A_{1}]-\ldots-[A_{q}].
\end{equation*}
Since $q>p$, without loss of generality 
$[A_{q}]\geq 0$ in $H_{1}({\mathbb 
R}^{3}-Y_{2};{\mathbb Z})$.  Hence
\begin{equation*}
-[A_{1}]-\ldots-[A_{q-1}]\geq p
\end{equation*}
in $H_{1}({\mathbb R}^{3}-Y_{2};{\mathbb Z})$.

In $H_{1}({\mathbb R}^{3}-Y_{1};{\mathbb Z})$ we have
\begin{equation*}
[X_{1}]=-[A_{1}]-\ldots-[A_{q}].
\end{equation*}
First we suppose that $[A_{q}]=0$ in $H_{1}({\mathbb 
R}^{3}-Y_{1};{\mathbb Z})$.  So 
$-[A_{1}]-\ldots-[A_{q-1}]=[X_{1}]\geq 1$ in 
$H_{1}({\mathbb R}^{3}-Y_{1};{\mathbb Z})$.  In this 
case, we let $L$ denote the simple closed curve 
obtained from $A_{1}\cup \ldots\cup A_{q-1}$ by 
omitting the edges $\overline {v_{j}w_{j}}$ for 
$j=2$, \dots, $q-1$.  We orient $L$ so that going
from $w_{1}$ to $v_{1}$ along the edge $\overline
{w_{1}v_{1}}$ is the  positive direction. Then
$L$ has at least
$2q$  vertices and $lk(L,Y_{1}) 
=lk(X_{1},Y_{1})\geq 1$ and $ \vert 
lk(L,Y_{2}) \vert \geq p$.  We are done by letting $Z=Y_1$ and
$W=Y_2$.

Now suppose that $[A_{q}]\mathbin{\not =}0$ in 
$H_{1}({\mathbb R}^{3}-Y_{1};{\mathbb Z})$.  Let 
$L$ denote the simple closed curve obtained from 
$X_{2}\cup A_{q}$ by omitting the edge $\overline 
{w_{q}w_{1}}$.  We orient $L$ so that going
from $v_{1}$ to $w_{1}$ along the edge $\overline
{v_{1}w_{1}}$ is the  positive direction. Then $L$ has at least $q+2$ 
vertices.  Also since $lk(X_{2},Y_{1})=0$ we have 
$ \vert lk(L,Y_{1}) \vert = \vert lk(A_{q},Y_{1}) 
\vert \geq 1$; Now by changing the orientation on $Y_1$, if
necessary, $lk(L,Y_{1}) \geq 1$.  Since
$[A_{q}]\geq  0$ in $H_{1}({\mathbb R}^{3}-Y_{2};{\mathbb Z})$ we 
have $lk(L,Y_{2})\geq 
lk(X_{2},Y_{2})=p$.  So we are done by letting $Z=Y_1$ and
$W=Y_2$.
\end{proof}

Now we will use Lemmas~\ref{lem1} and \ref{lem2},
together with an  inductive argument, to prove
Theorem~2.

\begin{thm2} Let $p\in {\mathbb N}$ be 
given, and let $n=p(15p-9)$.  Then every 
embedding of $K_{n}$ in ${\mathbb R}^{3}$ contains a 
$2$-component link $L_{p}\cup J_{p}$ such that 
for some orientation of $L_{p}\cup J_{p}$, we 
have $lk(L_{p},J_{p})\geq p$. 
\end{thm2}  

\begin{proof}   Suppose that for every oriented 
link $L_{p}\cup J_{p}$ contained in $K_{n}$ we 
have $lk(L_{p},J_{p})<p$.

	Let $G_{1}$, \dots, $G_{p}$ be $p$ disjoint copies 
of $K_{15p-9}$ which are contained in $K_{n}$.  
For each $i=1$, \dots, $p$, let $H_{i}$ be a 
subgraph of $G_{i}$ containing all $15p-9$ 
vertices of $G_{i}$ such that, as a topological 
space, $H_{i}$ is homeomorphic to $K_{6}$, yet 
between every pair of vertices which have valence 
$5$ in $H_{i}$ there is a path in $H_{i}$ 
containing $p-1$ vertices which are each of 
valence $2$ in $H_{i}$.  Since $K_{6}$ contains 
$15$ edges and $15(p-1)+6=15p-9$, there is such a 
subgraph $H_{i}$ in $K_{15p-9}$.  Now by 
\cite{CG}, each $H_{i}$ contains a link 
$X_{i}\cup Y_{i}$ such that with some orientation 
we have $lk(X_{i},Y_{i})\geq 1$ and $X_{i}$ 
and $Y_{i}$ each contain $3$ vertices with 
valence $5$ and $3(p-1)$ vertices with valence 
$2$ in $H_{i}$.  Thus $X_{i}$ and $Y_{i}$ each 
contain a total of $3p$ vertices.    

	We will prove by induction that for every 
$m=1$, \dots, $p$, the $K_{m(15p-9)}$, which has all 
of its vertices in $\bigcup _{i=1}^{m}G_{i}$, 
contains a link $L_{m}\cup J_{m}$ such that, with 
some orientation, $lk(L_{m},J_{m})=p_{m}\geq 
m$ and $L_{m}$ and $J_{m}$ each have at least 
$3p$ vertices.  We saw above that this is true 
for $m=1$.  Assume that it's true for $m$.  Thus 
$lk(L_{m},J_{m})=p_{m}$, and by our initial assumption 
$p_{m}<p$.  Also $G_{m+1}$ is disjoint from 
$K_{m(15p-9)}$ and contains a pair of simple 
closed curves $X_{m+1}$ and $Y_{m+1}$ each with 
$3p$ vertices such that 
$lk(X_{m+1},Y_{m+1})\geq 1$.  Thus the 
disjoint simple closed curves $L_{m}$, $J_{m}$, 
$X_{m+1}$, and $Y_{m+1}$ each contain at least 
$3p$ vertices and $3p>p_{m}$ since $p>p_{m}$.  
Thus by Lemma 2, the $K_{(m+1)(15p-9)}$, which 
has all of its vertices in $\bigcup 
_{i=1}^{m+1}G_{i}$, contains simple closed curves 
$L_{m+1}$, $Z_{m+1}$ and $W_{m+1}$ each with at 
least $3p$  vertices such that, for some 
orientation of $L_{m+1}\cup Z_{m+1}\cup W_{m+1}$, 
we have $lk(L_{m+1},Z_{m+1})=q_{m}\geq 1$ 
and $lk(L_{m+1},W_{m+1})=r_{m}\geq 
p_{m}\geq m$.  By assumption $q_{m}<p$ and 
$r_{m}<p$, so $q_{m}+r_{m}<3p$.  Thus we can 
apply Lemma 1 to obtain a simple closed curve 
$J_{m+1}$ which is contained in 
$K_{(m+1)(15p-9)}$ such that, with some 
orientation, $lk(L_{m+1},J_{m+1})\geq 
q_{m}+r_{m}\geq 1+m$ and $J_{m+1}$ has at 
least $6p$ vertices.  Thus for all $m=1$, \dots, $p$, 
the $K_{m(15p-9)}$, which has all of its vertices 
in $\bigcup _{i=1}^{m}G_{i}$, contains a link 
$L_{m}\cup J_{m}$ such that, with some 
orientation, $lk(L_{m},J_{m})\geq m$ and 
$L_{m}$ and $J_{m}$ each have at least $3p$ 
vertices.  It follows that $K_{p(15p-9)}$ 
contains a link $L_{p}\cup J_{p}$ such that, with 
some orientation, $lk(L_{p},J_{p})\geq p$.  
However this contradicts our assumption that for 
every oriented link $L_{p}\cup J_{p}$ contained 
in $K_{n}$ we have $lk(L_{p},J_{p})<p$.  Hence we 
must, in fact, have a $2$-component link 
$L_{p}\cup J_{p}$ in $K_{n}$ such that for some 
orientation of the link we have 
$lk(L_{p},J_{p})\geq p$. 
\end{proof} 
 
\begin{figure}[ht!]
\label{fig2}
\centerline{\small
%\ShowGrid 
\SetLabels 
\E(0.12*-.15){Positive crossing}\\
\E(.9*-.15){Negative crossing}\\
\endSetLabels 
\AffixLabels{{\includegraphics{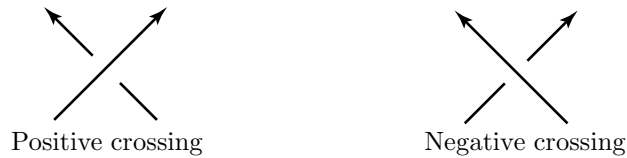}}}}
\vspace{3mm}
\caption{Positive and negative crossings}
\end{figure}

We now observe as follows that for any given $p\in {\mathbb N}$, 
it is not possible to have an $n\in {\mathbb N}$, such that 
every embedding of $K_{n}$ contains a 2-component 
link with linking number equal to precisely $p$.  
 First recall that the linking number of any 
oriented 2-component link is equal to one half of 
the sum of $+1$ for every positive crossing and 
$-1$ for every negative crossing.  We illustrate 
these two types of crossings in
Figure~\ref{fig2}.      

Let $n\in {\mathbb N}$ be fixed and let $K_{n}$ be 
embedded in ${\mathbb R}^{3}$.  We shall find a 
(possibly different) embedding of $K_{n}$ in 
${\mathbb R}^{3}$ which contains no $2$-component 
link $L_{1}\cup L_{2}$ with $lk(L_{1},L_{2})=p$.  
Let $q=max\lbrace  \vert lk(L_{i},L_{j}) \vert 
\rbrace $ taken over all disjoint pairs of simple 
closed curves $L_{i}$ and $L_{j}$, contained in 
$K_{n}$.  Suppose there is some link $L_{1}\cup 
L_{2}$ in $K_{n}$ such that for some orientation, 
$lk(L_{1},L_{2})=p$.  Let $e$ be an edge in 
$L_{1}$ and let $f$ be an edge in $L_{2}$.  We 
create a new embedding of $K_{n}$ by adding 
$2p+2q +2$ half twists between the edges $e$ and 
$f$ without changing anything else about the 
embedding.  We call this new embedding $K_{n}'$ 
and this new pair of edges $e'$ and $f'$.  Now 
for any pair of disjoint oriented simple closed 
curves $L_{i}$ and $L_{j}$ in $K_{n}$, there is a 
corresponding pair of disjoint oriented simple 
closed curves $L_{i}'$ and $L_{j}'$ in $K_{n}'$.  
Furthermore, $L_{i}'$ and $L_{j}'$ contain $e'$ 
and $f'$ respectively if and only if $L_{i}$ and 
$L_{j}$ contain $e$ and $f$ respectively.  Now 
suppose that one of $L_{i}'$ and $L_{j}'$ 
contains $e'$ and the other contains $f'$.  Then 
$lk(L_{i}',L_{j}')=lk(L_{i},L_{j})\pm (p+q+1)$.  
Since $lk(L_{i},L_{j})\geq -q$ we know that 
$lk(L_{i},L_{j})+(p+q+1)\geq p+1>0$, so $ 
\vert lk(L_{i},L_{j})+(p+q+1) \vert \geq 
p+1$.  Since $lk(L_{i},L_{j})\leq q$  we 
know that $lk(L_{i},L_{j})-(p+q+1)\leq 
-p-1<0$, so $ \vert lk(L_{i},L_{j})-(p+q+1) \vert 
\geq p+1$.  Thus $ \vert lk(L_{i}',L_{j}') 
\vert = \vert lk(L_{i},L_{j})\pm (p+q+1) \vert 
\geq p+1$.  

On the other hand, for any pair $L_{i}'$ and 
$L_{j}'$ which do not contain $e'$ and $f'$, then 
$lk(L_{i}',L_{j}')=lk(L_{i},L_{j})$.  It follows 
that $K_{n}'$ contains fewer links which have 
linking number equal to $p$ than $K_{n}$ did.  If 
we repeat this process enough times, we will eventually have 
an embedding of $K_{n}$ which contains no links 
that have linking number equal to $p$.   

\vspace{-1.5mm}
\section{Intrinsic
Knotting}  
\vspace{-1.5mm}

\begin{figure}[ht!]
\label{fig3}
\vspace{-2mm}
\centerline{\small
%\ShowGrid 
\SetLabels 
\E(0.5*-.05){$e_1$}\\
\E(0.5*.2){$e_0$}\\
\E(0.5*.8){$c_0$}\\
\E(0.5*1.03){$c_1$}\\
\E(-0.05*.5){$b_1$}\\
\E(.2*.5){$b_0$}\\
\E(.8*.5){$d_0$}\\
\E(1.05*.5){$d_1$}\\
\endSetLabels 
\AffixLabels{{\includegraphics{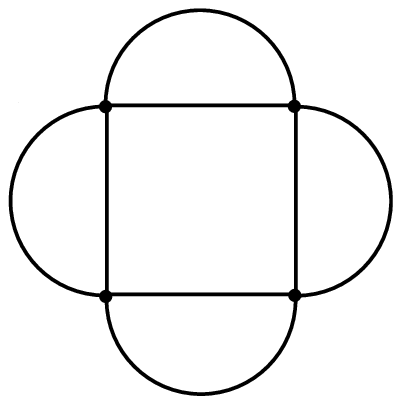}}}}
\vspace{2mm}
\caption{The pseudo-graph $G_4$}
\end{figure}

In order to consider intrinsic knotting of 
complete graphs, we shall make use of the 
pseudo-graph $D_{4}$ (illustrated in
Figure~\ref{fig3}) in  a similar manner to that of
\cite{TY} and 
\cite{Fo}.

 A \emph{Hamiltonian cycle} in a graph or 
pseudo-graph $G$ is a simple closed curve in $G$ 
which contains every vertex of $G$.  Each 
Hamiltonian cycle of $D_{4}$ has the form 
$b_{i}c_{j}d_{k}e_{l}$ where $i$, $j$, $k$, $l\in 
\lbrace 0,1\rbrace $.  For each  such cycle we 
define 
\begin{equation*}
\varepsilon (b_{i}c_{j}d_{k}e_{l})=
\begin{cases} 
+1 & \text{if } i+j+k+l \text 
{ is even} \\
-1& \text{if } i+j+k+l \text 
{ is odd}
\end{cases}
\end{equation*}
Let $\Gamma $ denote the set of all Hamiltonian 
cycles in $D_{4}$.  Recall that if $\gamma $ is a 
simple closed curve in ${\mathbb R}^{3}$, then 
$a_{2}(\gamma )$ denotes the second coefficient 
of the Conway polynomial of $\gamma $.  For any 
embedding of $D_{4}$ in ${\mathbb R}^{3}$ we define 
\begin{equation*}
\lambda = \vert \sum _{\gamma \in
\Gamma  }\varepsilon (\gamma )a_{2}(\gamma ) \vert 
\end{equation*}
Let $B$, $C$, $D$, and $E$ denote the simple 
closed curves $b_{0}\cup b_{1}$, $c_{0}\cup 
c_{1}$, $d_{0}\cup d_{1}$, and $e_{0}\cup e_{1}$ 
respectively.  It is shown in \cite{TY} that, for 
a given embedding of $D_{4}$ in ${\mathbb R}^{3}$, 
we have the equation
\begin{equation}
\label{eqn2}
\lambda = \vert lk(E,C)lk(B,D)
\vert      
\end{equation}

Equation (2) enables us to relate knotting and 
linking in an embedded graph.  We use this 
equation to prove Theorem 3.

\begin{thm3}   Let $m\in {\mathbb N}$ be 
given, let $p$ denote an integer such that 
$p\geq 4\sqrt {m}$, and let $n=p(15p-9)$.  
Then every embedding of $K_{2n}$ in ${\mathbb 
R}^{3}$ contains a knot $Q$ with $ \vert a_{2}(Q) 
\vert \geq m$. 
\end{thm3}    

\begin{proof}     Let $K_{2n}$ be embedded in 
${\mathbb R}^{3}$.  Let $G_{1}$ and $G_{2}$ denote 
two disjoint $K_{n}$'s which are contained in 
$K_{2n}$.  By Theorem 2, each $G_{i}$ contains a 
link $X_{i}\cup Y_{i}$ such that, for some 
orientation of $X_{i}\cup Y_{i}$, we have 
$lk(X_{i},Y_{i})\geq p\geq 4\sqrt {m}$. 
 For $i=1$ and $i=2$, let $u_{i}$ and $v_{i}$ be 
distinct vertices on $X_{i}$, and let $z_{i}$ and 
$w_{i}$ be distinct vertices on $Y_{i}$.  
Consider the embedded subgraph $G$ of $K_{2n}$ consisting 
of $X_{1}\cup Y_{1}\cup X_{2}\cup Y_{2}$ together 
with the edges $\overline {v_{1}u_{2}}$, 
$\overline {w_{1}v_{2}}$, $\overline 
{z_{1}w_{2}}$, and $\overline {u_{1}z_{2}}$.  
Figure~4 illustrates the abstract graph
$G$,  with all of the other vertices omitted.

\begin{figure}[ht!]
\label{fig4}
\vspace{3mm}
\centerline{\small
%\ShowGrid 
\SetLabels 
\E(-0.05*.5){$X_1$}\\
\E(.11*.3){$u_1$}\\
\E(.11*.67){$v_1$}\\
\E(.32*.13){$z_2$}\\
\E(.32*.87){$u_2$}\\
\E(.5*-.05){$Y_2$}\\
\E(.5*1.05){$X_2$}\\
\E(.69*.13){$w_2$}\\
\E(.67*.87){$v_2$}\\
\E(.9*.3){$z_1$}\\
\E(.9*.67){$w_1$}\\
\E(1.05*.5){$Y_1$}\\
\endSetLabels 
\AffixLabels{{\includegraphics{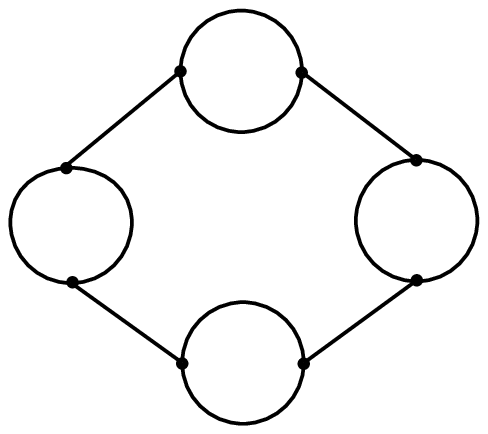}}}}
\vspace{5mm}
\caption{The graph $G$}
\end{figure}  

 From the embedded graph $G$, we obtain an embedding of 
the pseudo-graph $D_{4}$ by omitting all vertices 
of valence $2$ and collapsing the edges 
$\overline {v_{1}u_{2}}$, $\overline 
{w_{1}v_{2}}$, $\overline {z_{1}w_{2}}$, and 
$\overline {u_{1}z_{2}}$.  Now by equation $(2)$ 
we have $\lambda = \vert 
lk(X_{1},Y_{1})lk(X_{2},Y_{2}) \vert \geq 
16m$.  Thus $\sum _{\gamma \in \Gamma } \vert 
a_{2}(\gamma ) \vert \geq  \vert \sum 
_{\gamma \in \Gamma }\varepsilon (\gamma 
)a_{2}(\gamma ) \vert \geq 16m$.  However 
$D_{4}$ has precisely $16$ Hamiltonian cycles.  
Thus there is some $\gamma _{0}\in \Gamma $ such 
that $ \vert a_{2}(\gamma _{0}) \vert \geq 
m$.  Now observe that there is a simple closed 
$Q$ in $K_{2n}$ such that $\gamma _{0}$ is 
obtained from $Q$ by collapsing precisely the 
edges $\overline {v_{1}u_{2}}$, $\overline 
{w_{1}v_{2}}$, $\overline {z_{1}w_{2}}$, and 
$\overline {u_{1}z_{2}}$.  Thus $Q$ has the same 
knot type as $\gamma _{0}$ and hence $ \vert 
a_{2}(Q) \vert \geq m$.
\end{proof}

We can see as follows that for a given $m\in 
{\mathbb N}$ there is no $n\in {\mathbb N}$ such that 
every embedding of $K_{n}$ contains a knot $Q$ 
with $ \vert a_{2}(Q) \vert =m$.  First observe 
that if $Q$ is any knot and $R$ is a right handed 
trefoil knot, then it follows from equation $(1)$ 
in the introduction that 
$a_{2}(Q\#R)=a_{2}(Q)+1$.  Let $m\in {\mathbb N}$ 
and $n\in {\mathbb N}$ be fixed, and let $K_{n}$ be 
embedded in ${\mathbb R}^{3}$.  Consider all simple 
closed curves $Q_{i}$ which are contained in 
$K_{n}$, and let $q=min\lbrace 
a_{2}(Q_{i}),0\rbrace \leq 0$.  We modify 
the embedding of $K_{n}$ by adding $m-q\geq 
m$ right handed trefoil knots to every edge 
within a ball which does not intersect the graph 
elsewhere.  We call this modified embedding 
$K_{n}'$.  Now for each simple closed curve 
$Q_{i}'$ in $K_{n}'$ there is a corresponding 
simple closed curve $Q_{i}$ in $K_{n}$.  Since 
every simple closed curve in $K_{n}'$ contains at 
least $3$ edges, each $Q_{i}'$ is a connected sum 
which has at least $3(m-q)$ more right-handed 
trefoil knots as summands than $Q_{i}$ had.  Thus 
for every simple closed curve $Q_{i}'$ in $K_{n}'$, we have
$a_{2}(Q_{i}')\geq  a_{2}(Q_{i})+3(m-q)\geq 3m-2q\geq 3m$.  Thus
$K_{n}'$ contains no  knot $Q'$ such that $ \vert a_{2}(Q') \vert =m$.

{\makeatletter
\@thebibliography@{FNP}\small\parskip0pt % 
plus2pt\relax
\makeatother

\bibitem [CG] {CG} \textbf{J. Conway}, \textbf{C.
McA Gordon},
\emph{Knots and links in spatial graphs},  J. of 
Graph Theory 7 (1983), 445--453.

\bibitem [FNP]{FNP} \textbf{E. Flapan}, \textbf{R.
Naimi}, \textbf{J.  Pommersheim},
\emph{Intrinsically triple linked  complete
graphs}, Topology and its  Applications
115 (2001), 239--246.

\bibitem [Fo]{Fo} \textbf{J. Foisy},
\emph{Intrinsically  knotted graphs}, J. of Graph
Theory  39 (2002), 178--187.

\bibitem [Ka]{Ka} \textbf{L.H. Kauffman}, \emph{The
Conway  polynomial}, Topology 20
(1981), 101--108.

\bibitem [Sh]{Sh} \textbf{M. Shimabara}, \emph{Knots
in  certain spatial graphs}, Tokyo J. Math.  
11 (1988), 405--413.

\bibitem [TY]{TY} \textbf{K. Taniyama}, \textbf{A.
Yasuhara}, \emph{Realization of knots and links in a
spatial  graph}, Topology and its
Applications  112 (2001), 87--109.

\endthebibliography}

\Addressesr

\end{document}

%% file: agtout.tex
\def\ifplaintex{\expandafter\ifx\csname documentclass\endcsname\relax}

\def\gtp{{\mathsurround=0pt\it $\cal G\mskip-2mu$eometry \&\ 
$\cal T\!\!$opology $\cal P\!$ublications}}  % GT publications

\def\Addressesr{\bigskip
{\small \parskip 0pt \leftskip 0pt \rightskip 0pt plus 1fil \def\\{\par}
\sl\theaddress\par
\medskip
\rm Email:\stdspace\tt\theemail\hfill\rm Received:\qua\receiveddate \par}}

\def\recd{{\small Received:\qua\receiveddate\ifx\reviseddate\relax
\else\qquad Revised:\qua\reviseddate\fi\par}} 

\def\lognumber#1{\def\thelognumber{#1}}
\def\volumenumber#1{\def\thevolumenumber{#1}}
\def\volumeyear#1{\def\thevolumeyear{#1}}
\def\papernumber#1{\def\thepapernumber{#1}}
\def\pagenumbers#1#2{\def\startpage{#1}\def\finishpage{#2}}
\def\published#1{\def\publishdate{#1}}

\def\received#1{\def\receiveddate{#1}}

\def\accepted#1{\def\accepteddate{#1}}

\long\def\asciiabstract#1{\long\def\theasciiabstract{#1}}

\let\\\par\let\thelognumber\relax\let\thevolumenumber\relax
\let\thepapernumber\relax\let\thevolumeyear\relax\let\startpage\relax
\let\finishpage\relax\let\publishdate\relax\let\receiveddate\relax
\let\reviseddate\relax\let\accepteddate\relax\let\theasciititle\relax
\let\theasciiauthors\relax
\let\theasciiabstract\relax

\let\theasciiemail\relax

\ifplaintex
\font\logobig=cmssbx10 scaled 3836
\font\logomed=cmssbx10 scaled 2557
\else
\font\logobig=cmssbx10 scaled 4200
\font\logomed=cmssbx10 scaled 2800
\fi

\long\def\makeagttitle{   %%% start of definition of \makeagttitle
\count0=\startpage
\agt\hfill      %   Journal title (top left) 
\hbox to 45truept{\vbox to 0pt{\vglue -13truept{\logomed A\kern -.37em{\logobig 
T}\kern -.38em G}\vss}\hss}
\break
{\small Volume \thevolumenumber\ (\thevolumeyear)
\startpage--\finishpage\nl
Published: \publishdate}

\vglue .25truein

{\parskip=0pt\leftskip 0pt plus
1fil\def\\{\par\smallskip}{\Large\bf\thetitle}\par\medskip} \vglue
0.05truein

{\parskip=0pt\leftskip 0pt plus 1fil\def\\{\par}{\sc\theauthors}
\par\medskip}%
 
\vglue 0.03truein 

{\small\leftskip 25truept\rightskip 25truept{\bf Abstract}\stdspace\theabstract

{\bf AMS Classification}\stdspace\theprimaryclass
\ifx\thesecondaryclass\relax\else; \thesecondaryclass\fi\par
{\bf Keywords}\stdspace \thekeywords\par}\vglue 7truept

}   %%%% end of definition of \makeagttitle

\ifplaintex
\font\phead=cmsl9 scaled 950
\font\pnum=cmbx10 scaled 913
\font\pfoot=cmsl9 scaled 950
\headline{\vbox to 0pt{\vskip -4.5mm\line{\small\phead\ifnum
\count0=\startpage ISSN 1472-2739 (on-line) 1472-2747 (printed)
\hfill {\pnum\folio}\else\ifodd\count0\def\\{ }% 
\ifx\theshorttitle\relax\thetitle\else\theshorttitle\fi\hfill{\pnum\folio}
\else\def\\{ and }{\pnum\folio}\hfill\ifx\theshortauthors\relax\theauthors
\else\theshortauthors\fi\fi\fi}\vss}}
\footline{\vbox to 0pt{\vglue 0mm\line{\small\pfoot\ifnum\count0=\startpage
\copyright\ \gtp\hfill\else
\agt, Volume \thevolumenumber\ (\thevolumeyear)\hfill\fi}\vss}}
\else
\font=cmsl9 scaled 1050
\font\lnum=cmbx10 
\font=cmsl9 scaled 1050
\makeatletter
\def\@oddhead{{\small\lhead\ifnum\count0=\startpage ISSN 1472-2739 
(on-line) 1472-2747 (printed)\hfill {\lnum\number\count0}\else\ifodd\count0
\def\\{ }\ifx\theshorttitle\relax \thetitle \else\theshorttitle\fi\hfill
{\lnum\number\count0}\else\def\\{ and }{\lnum\number\count0}
\hfill\ifx\theshortauthors\relax 
\theauthors\else\theshortauthors\fi\fi\fi}}\def\@evenhead{\@oddhead}
\def\@oddfoot{\small\lfoot\ifnum\count0=\startpage\copyright\ \gtp\hfill\else
\agt, Volume \thevolumenumber\ (\thevolumeyear)\hfill\fi}
\def\@evenfoot{\@oddfoot}
\makeatother
\fi
\let\maketitlepage\makeagttitle
\let\makeshorttitle\maketitlepage
\let\maketitle\maketitlepage

\newwrite\gtoutfile
\long\gdef\makeheadfile{  %%% start of definition of \makeheadfile
{\def\\{, }\def\s{ }
\immediate\openout\gtoutfile head.xxx
\immediate\write\gtoutfile{To: math@arxiv.org}
\immediate\write\gtoutfile{Subject: put OR rep NNNNN:ppppp}
\immediate\write\gtoutfile{--text follows this line--}
\immediate\write\gtoutfile{Proxy-for: \ifx\theasciiauthors\relax
\theauthors\else\theasciiauthors\fi\s<\ifx\theasciiemail\relax\theemail\else\theasciiemail\fi>}
\immediate\write\gtoutfile{\noexpand\\}
\immediate\write\gtoutfile{Authors: \ifx\theasciiauthors\relax
\theauthors\else\theasciiauthors\fi}
{\def\\{ }\immediate\write\gtoutfile{Title: \ifx\theasciititle\relax
\thetitle\else\theasciititle\fi}}
\immediate\write\gtoutfile{Subj-class: GT or SG, GR etc}
\immediate\write\gtoutfile{MSC-class: \theprimaryclass\ifx\thesecondaryclass\relax\else, \thesecondaryclass\fi}
\immediate\write\gtoutfile{Journal-ref: Algebr. Geom. Topol. \thevolumenumber\s
(\thevolumeyear) \startpage-\finishpage}
\immediate\write\gtoutfile{Comments: Published by Algebraic and
Geometric Topology at}
\immediate\write\gtoutfile{\s\s\s  http://www.maths.warwick.ac.uk/agt/AGTVol\thevolumenumber/agt-\thevolumenumber-\thepapernumber.abs.html}
\immediate\write\gtoutfile{\noexpand\\}
\immediate\write\gtoutfile{}
\ifx\theasciiabstract\relax
\immediate\write\gtoutfile{\theabstract}\else
\immediate\write\gtoutfile{\theasciiabstract}\fi
\immediate\write\gtoutfile{}
\immediate\write\gtoutfile{\noexpand\\}
\immediate\write\gtoutfile{}
\immediate\closeout\gtoutfile}}  %%% end of definition of \makeheadfile

\def\maketitlepage{\makeagttitle\makeheadfile}
\let\makeshorttitle\maketitlepage
\let\maketitle\maketitlepage

\def\ifplaintex{\expandafter\ifx\csname documentclass\endcsname\relax}

\def\gtp{{\mathsurround=0pt\it $\cal G\mskip-2mu$eometry \&\ 
$\cal T\!\!$opology $\cal P\!$ublications}}  % GT publications

\def\Addressesr{\bigskip
{\small \parskip 0pt \leftskip 0pt \rightskip 0pt plus 1fil \def\\{\par}
\sl\theaddress\par
\medskip
\rm Email:\stdspace\tt\theemail\hfill\rm Received:\qua\receiveddate \par}}

\def\recd{{\small Received:\qua\receiveddate\ifx\reviseddate\relax
\else\qquad Revised:\qua\reviseddate\fi\par}} 

\def\lognumber#1{\def\thelognumber{#1}}
\def\volumenumber#1{\def\thevolumenumber{#1}}
\def\volumeyear#1{\def\thevolumeyear{#1}}
\def\papernumber#1{\def\thepapernumber{#1}}
\def\pagenumbers#1#2{\def\startpage{#1}\def\finishpage{#2}}
\def\published#1{\def\publishdate{#1}}

\def\received#1{\def\receiveddate{#1}}

\def\accepted#1{\def\accepteddate{#1}}

\long\def\asciiabstract#1{\long\def\theasciiabstract{#1}}

\let\\\par\let\thelognumber\relax\let\thevolumenumber\relax
\let\thepapernumber\relax\let\thevolumeyear\relax\let\startpage\relax
\let\finishpage\relax\let\publishdate\relax\let\receiveddate\relax
\let\reviseddate\relax\let\accepteddate\relax\let\theasciititle\relax
\let\theasciiauthors\relax
\let\theasciiabstract\relax

\let\theasciiemail\relax

\ifplaintex
\font\logobig=cmssbx10 scaled 3836
\font\logomed=cmssbx10 scaled 2557
\else
\font\logobig=cmssbx10 scaled 4200
\font\logomed=cmssbx10 scaled 2800
\fi

\long\def\makeagttitle{   %%% start of definition of \makeagttitle
\count0=\startpage
\agt\hfill      %   Journal title (top left) 
\hbox to 45truept{\vbox to 0pt{\vglue -13truept{\logomed A\kern -.37em{\logobig 
T}\kern -.38em G}\vss}\hss}
\break
{\small Volume \thevolumenumber\ (\thevolumeyear)
\startpage--\finishpage\nl
Published: \publishdate}

\vglue .25truein

{\parskip=0pt\leftskip 0pt plus
1fil\def\\{\par\smallskip}{\Large\bf\thetitle}\par\medskip} \vglue
0.05truein

{\parskip=0pt\leftskip 0pt plus 1fil\def\\{\par}{\sc\theauthors}
\par\medskip}%
 
\vglue 0.03truein 

{\small\leftskip 25truept\rightskip 25truept{\bf Abstract}\stdspace\theabstract

{\bf AMS Classification}\stdspace\theprimaryclass
\ifx\thesecondaryclass\relax\else; \thesecondaryclass\fi\par
{\bf Keywords}\stdspace \thekeywords\par}\vglue 7truept

}   %%%% end of definition of \makeagttitle

\ifplaintex
\font\phead=cmsl9 scaled 950
\font\pnum=cmbx10 scaled 913
\font\pfoot=cmsl9 scaled 950
\headline{\vbox to 0pt{\vskip -4.5mm\line{\small\phead\ifnum
\count0=\startpage ISSN 1472-2739 (on-line) 1472-2747 (printed)
\hfill {\pnum\folio}\else\ifodd\count0\def\\{ }% 
\ifx\theshorttitle\relax\thetitle\else\theshorttitle\fi\hfill{\pnum\folio}
\else\def\\{ and }{\pnum\folio}\hfill\ifx\theshortauthors\relax\theauthors
\else\theshortauthors\fi\fi\fi}\vss}}
\footline{\vbox to 0pt{\vglue 0mm\line{\small\pfoot\ifnum\count0=\startpage
\copyright\ \gtp\hfill\else
\agt, Volume \thevolumenumber\ (\thevolumeyear)\hfill\fi}\vss}}
\else
\font=cmsl9 scaled 1050
\font\lnum=cmbx10 
\font=cmsl9 scaled 1050
\makeatletter
\def\@oddhead{{\small\lhead\ifnum\count0=\startpage ISSN 1472-2739 
(on-line) 1472-2747 (printed)\hfill {\lnum\number\count0}\else\ifodd\count0
\def\\{ }\ifx\theshorttitle\relax \thetitle \else\theshorttitle\fi\hfill
{\lnum\number\count0}\else\def\\{ and }{\lnum\number\count0}
\hfill\ifx\theshortauthors\relax 
\theauthors\else\theshortauthors\fi\fi\fi}}\def\@evenhead{\@oddhead}
\def\@oddfoot{\small\lfoot\ifnum\count0=\startpage\copyright\ \gtp\hfill\else
\agt, Volume \thevolumenumber\ (\thevolumeyear)\hfill\fi}
\def\@evenfoot{\@oddfoot}
\makeatother
\fi
\let\maketitlepage\makeagttitle
\let\makeshorttitle\maketitlepage
\let\maketitle\maketitlepage

\newwrite\gtoutfile
\long\gdef\makeheadfile{  %%% start of definition of \makeheadfile
{\def\\{, }\def\s{ }
\immediate\openout\gtoutfile head.xxx
\immediate\write\gtoutfile{To: math@arxiv.org}
\immediate\write\gtoutfile{Subject: put OR rep NNNNN:ppppp}
\immediate\write\gtoutfile{--text follows this line--}
\immediate\write\gtoutfile{Proxy-for: \ifx\theasciiauthors\relax
\theauthors\else\theasciiauthors\fi\s<\ifx\theasciiemail\relax\theemail\else\theasciiemail\fi>}
\immediate\write\gtoutfile{\noexpand\\}
\immediate\write\gtoutfile{Authors: \ifx\theasciiauthors\relax
\theauthors\else\theasciiauthors\fi}
{\def\\{ }\immediate\write\gtoutfile{Title: \ifx\theasciititle\relax
\thetitle\else\theasciititle\fi}}
\immediate\write\gtoutfile{Subj-class: GT or SG, GR etc}
\immediate\write\gtoutfile{MSC-class: \theprimaryclass\ifx\thesecondaryclass\relax\else, \thesecondaryclass\fi}
\immediate\write\gtoutfile{Journal-ref: Algebr. Geom. Topol. \thevolumenumber\s
(\thevolumeyear) \startpage-\finishpage}
\immediate\write\gtoutfile{Comments: Published by Algebraic and
Geometric Topology at}
\immediate\write\gtoutfile{\s\s\s  http://www.maths.warwick.ac.uk/agt/AGTVol\thevolumenumber/agt-\thevolumenumber-\thepapernumber.abs.html}
\immediate\write\gtoutfile{\noexpand\\}
\immediate\write\gtoutfile{}
\ifx\theasciiabstract\relax
\immediate\write\gtoutfile{\theabstract}\else
\immediate\write\gtoutfile{\theasciiabstract}\fi
\immediate\write\gtoutfile{}
\immediate\write\gtoutfile{\noexpand\\}
\immediate\write\gtoutfile{}
\immediate\closeout\gtoutfile}}  %%% end of definition of \makeheadfile

\def\maketitlepage{\makeagttitle\makeheadfile}
\let\makeshorttitle\maketitlepage
\let\maketitle\maketitlepage

\def\ifplaintex{\expandafter\ifx\csname documentclass\endcsname\relax}

\def\gtp{{\mathsurround=0pt\it $\cal G\mskip-2mu$eometry \&\ 
$\cal T\!\!$opology $\cal P\!$ublications}}  % GT publications

\def\Addressesr{\bigskip
{\small \parskip 0pt \leftskip 0pt \rightskip 0pt plus 1fil \def\\{\par}
\sl\theaddress\par
\medskip
\rm Email:\stdspace\tt\theemail\hfill\rm Received:\qua\receiveddate \par}}

\def\recd{{\small Received:\qua\receiveddate\ifx\reviseddate\relax
\else\qquad Revised:\qua\reviseddate\fi\par}} 

\def\lognumber#1{\def\thelognumber{#1}}
\def\volumenumber#1{\def\thevolumenumber{#1}}
\def\volumeyear#1{\def\thevolumeyear{#1}}
\def\papernumber#1{\def\thepapernumber{#1}}
\def\pagenumbers#1#2{\def\startpage{#1}\def\finishpage{#2}}
\def\published#1{\def\publishdate{#1}}

\def\received#1{\def\receiveddate{#1}}

\def\accepted#1{\def\accepteddate{#1}}

\long\def\asciiabstract#1{\long\def\theasciiabstract{#1}}

\let\\\par\let\thelognumber\relax\let\thevolumenumber\relax
\let\thepapernumber\relax\let\thevolumeyear\relax\let\startpage\relax
\let\finishpage\relax\let\publishdate\relax\let\receiveddate\relax
\let\reviseddate\relax\let\accepteddate\relax\let\theasciititle\relax
\let\theasciiauthors\relax
\let\theasciiabstract\relax

\let\theasciiemail\relax

\ifplaintex
\font\logobig=cmssbx10 scaled 3836
\font\logomed=cmssbx10 scaled 2557
\else
\font\logobig=cmssbx10 scaled 4200
\font\logomed=cmssbx10 scaled 2800
\fi

\long\def\makeagttitle{   %%% start of definition of \makeagttitle
\count0=\startpage
\agt\hfill      %   Journal title (top left) 
\hbox to 45truept{\vbox to 0pt{\vglue -13truept{\logomed A\kern -.37em{\logobig 
T}\kern -.38em G}\vss}\hss}
\break
{\small Volume \thevolumenumber\ (\thevolumeyear)
\startpage--\finishpage\nl
Published: \publishdate}

\vglue .25truein

{\parskip=0pt\leftskip 0pt plus
1fil\def\\{\par\smallskip}{\Large\bf\thetitle}\par\medskip} \vglue
0.05truein

{\parskip=0pt\leftskip 0pt plus 1fil\def\\{\par}{\sc\theauthors}
\par\medskip}%
 
\vglue 0.03truein 

{\small\leftskip 25truept\rightskip 25truept{\bf Abstract}\stdspace\theabstract

{\bf AMS Classification}\stdspace\theprimaryclass
\ifx\thesecondaryclass\relax\else; \thesecondaryclass\fi\par
{\bf Keywords}\stdspace \thekeywords\par}\vglue 7truept

}   %%%% end of definition of \makeagttitle

\ifplaintex
\font\phead=cmsl9 scaled 950
\font\pnum=cmbx10 scaled 913
\font\pfoot=cmsl9 scaled 950
\headline{\vbox to 0pt{\vskip -4.5mm\line{\small\phead\ifnum
\count0=\startpage ISSN 1472-2739 (on-line) 1472-2747 (printed)
\hfill {\pnum\folio}\else\ifodd\count0\def\\{ }% 
\ifx\theshorttitle\relax\thetitle\else\theshorttitle\fi\hfill{\pnum\folio}
\else\def\\{ and }{\pnum\folio}\hfill\ifx\theshortauthors\relax\theauthors
\else\theshortauthors\fi\fi\fi}\vss}}
\footline{\vbox to 0pt{\vglue 0mm\line{\small\pfoot\ifnum\count0=\startpage
\copyright\ \gtp\hfill\else
\agt, Volume \thevolumenumber\ (\thevolumeyear)\hfill\fi}\vss}}
\else
\font=cmsl9 scaled 1050
\font\lnum=cmbx10 
\font=cmsl9 scaled 1050
\makeatletter
\def\@oddhead{{\small\lhead\ifnum\count0=\startpage ISSN 1472-2739 
(on-line) 1472-2747 (printed)\hfill {\lnum\number\count0}\else\ifodd\count0
\def\\{ }\ifx\theshorttitle\relax \thetitle \else\theshorttitle\fi\hfill
{\lnum\number\count0}\else\def\\{ and }{\lnum\number\count0}
\hfill\ifx\theshortauthors\relax 
\theauthors\else\theshortauthors\fi\fi\fi}}\def\@evenhead{\@oddhead}
\def\@oddfoot{\small\lfoot\ifnum\count0=\startpage\copyright\ \gtp\hfill\else
\agt, Volume \thevolumenumber\ (\thevolumeyear)\hfill\fi}
\def\@evenfoot{\@oddfoot}
\makeatother
\fi
\let\maketitlepage\makeagttitle
\let\makeshorttitle\maketitlepage
\let\maketitle\maketitlepage

\newwrite\gtoutfile
\long\gdef\makeheadfile{  %%% start of definition of \makeheadfile
{\def\\{, }\def\s{ }
\immediate\openout\gtoutfile head.xxx
\immediate\write\gtoutfile{To: math@arxiv.org}
\immediate\write\gtoutfile{Subject: put OR rep NNNNN:ppppp}
\immediate\write\gtoutfile{--text follows this line--}
\immediate\write\gtoutfile{Proxy-for: \ifx\theasciiauthors\relax
\theauthors\else\theasciiauthors\fi\s<\ifx\theasciiemail\relax\theemail\else\theasciiemail\fi>}
\immediate\write\gtoutfile{\noexpand\\}
\immediate\write\gtoutfile{Authors: \ifx\theasciiauthors\relax
\theauthors\else\theasciiauthors\fi}
{\def\\{ }\immediate\write\gtoutfile{Title: \ifx\theasciititle\relax
\thetitle\else\theasciititle\fi}}
\immediate\write\gtoutfile{Subj-class: GT or SG, GR etc}
\immediate\write\gtoutfile{MSC-class: \theprimaryclass\ifx\thesecondaryclass\relax\else, \thesecondaryclass\fi}
\immediate\write\gtoutfile{Journal-ref: Algebr. Geom. Topol. \thevolumenumber\s
(\thevolumeyear) \startpage-\finishpage}
\immediate\write\gtoutfile{Comments: Published by Algebraic and
Geometric Topology at}
\immediate\write\gtoutfile{\s\s\s  http://www.maths.warwick.ac.uk/agt/AGTVol\thevolumenumber/agt-\thevolumenumber-\thepapernumber.abs.html}
\immediate\write\gtoutfile{\noexpand\\}
\immediate\write\gtoutfile{}
\ifx\theasciiabstract\relax
\immediate\write\gtoutfile{\theabstract}\else
\immediate\write\gtoutfile{\theasciiabstract}\fi
\immediate\write\gtoutfile{}
\immediate\write\gtoutfile{\noexpand\\}
\immediate\write\gtoutfile{}
\immediate\closeout\gtoutfile}}  %%% end of definition of \makeheadfile

\def\maketitlepage{\makeagttitle\makeheadfile}
\let\makeshorttitle\maketitlepage
\let\maketitle\maketitlepage

\def\ifplaintex{\expandafter\ifx\csname documentclass\endcsname\relax}

\def\gtp{{\mathsurround=0pt\it $\cal G\mskip-2mu$eometry \&\ 
$\cal T\!\!$opology $\cal P\!$ublications}}  % GT publications

\def\Addressesr{\bigskip
{\small \parskip 0pt \leftskip 0pt \rightskip 0pt plus 1fil \def\\{\par}
\sl\theaddress\par
\medskip
\rm Email:\stdspace\tt\theemail\hfill\rm Received:\qua\receiveddate \par}}

\def\recd{{\small Received:\qua\receiveddate\ifx\reviseddate\relax
\else\qquad Revised:\qua\reviseddate\fi\par}} 

\def\lognumber#1{\def\thelognumber{#1}}
\def\volumenumber#1{\def\thevolumenumber{#1}}
\def\volumeyear#1{\def\thevolumeyear{#1}}
\def\papernumber#1{\def\thepapernumber{#1}}
\def\pagenumbers#1#2{\def\startpage{#1}\def\finishpage{#2}}
\def\published#1{\def\publishdate{#1}}

\def\received#1{\def\receiveddate{#1}}

\def\accepted#1{\def\accepteddate{#1}}

\long\def\asciiabstract#1{\long\def\theasciiabstract{#1}}

\let\\\par\let\thelognumber\relax\let\thevolumenumber\relax
\let\thepapernumber\relax\let\thevolumeyear\relax\let\startpage\relax
\let\finishpage\relax\let\publishdate\relax\let\receiveddate\relax
\let\reviseddate\relax\let\accepteddate\relax\let\theasciititle\relax
\let\theasciiauthors\relax
\let\theasciiabstract\relax

\let\theasciiemail\relax

\ifplaintex
\font\logobig=cmssbx10 scaled 3836
\font\logomed=cmssbx10 scaled 2557
\else
\font\logobig=cmssbx10 scaled 4200
\font\logomed=cmssbx10 scaled 2800
\fi

\long\def\makeagttitle{   %%% start of definition of \makeagttitle
\count0=\startpage
\agt\hfill      %   Journal title (top left) 
\hbox to 45truept{\vbox to 0pt{\vglue -13truept{\logomed A\kern -.37em{\logobig 
T}\kern -.38em G}\vss}\hss}
\break
{\small Volume \thevolumenumber\ (\thevolumeyear)
\startpage--\finishpage\nl
Published: \publishdate}

\vglue .25truein

{\parskip=0pt\leftskip 0pt plus
1fil\def\\{\par\smallskip}{\Large\bf\thetitle}\par\medskip} \vglue
0.05truein

{\parskip=0pt\leftskip 0pt plus 1fil\def\\{\par}{\sc\theauthors}
\par\medskip}%
 
\vglue 0.03truein 

{\small\leftskip 25truept\rightskip 25truept{\bf Abstract}\stdspace\theabstract

{\bf AMS Classification}\stdspace\theprimaryclass
\ifx\thesecondaryclass\relax\else; \thesecondaryclass\fi\par
{\bf Keywords}\stdspace \thekeywords\par}\vglue 7truept

}   %%%% end of definition of \makeagttitle

\ifplaintex
\font\phead=cmsl9 scaled 950
\font\pnum=cmbx10 scaled 913
\font\pfoot=cmsl9 scaled 950
\headline{\vbox to 0pt{\vskip -4.5mm\line{\small\phead\ifnum
\count0=\startpage ISSN 1472-2739 (on-line) 1472-2747 (printed)
\hfill {\pnum\folio}\else\ifodd\count0\def\\{ }% 
\ifx\theshorttitle\relax\thetitle\else\theshorttitle\fi\hfill{\pnum\folio}
\else\def\\{ and }{\pnum\folio}\hfill\ifx\theshortauthors\relax\theauthors
\else\theshortauthors\fi\fi\fi}\vss}}
\footline{\vbox to 0pt{\vglue 0mm\line{\small\pfoot\ifnum\count0=\startpage
\copyright\ \gtp\hfill\else
\agt, Volume \thevolumenumber\ (\thevolumeyear)\hfill\fi}\vss}}
\else
\font=cmsl9 scaled 1050
\font\lnum=cmbx10 
\font=cmsl9 scaled 1050
\makeatletter
\def\@oddhead{{\small\lhead\ifnum\count0=\startpage ISSN 1472-2739 
(on-line) 1472-2747 (printed)\hfill {\lnum\number\count0}\else\ifodd\count0
\def\\{ }\ifx\theshorttitle\relax \thetitle \else\theshorttitle\fi\hfill
{\lnum\number\count0}\else\def\\{ and }{\lnum\number\count0}
\hfill\ifx\theshortauthors\relax 
\theauthors\else\theshortauthors\fi\fi\fi}}\def\@evenhead{\@oddhead}
\def\@oddfoot{\small\lfoot\ifnum\count0=\startpage\copyright\ \gtp\hfill\else
\agt, Volume \thevolumenumber\ (\thevolumeyear)\hfill\fi}
\def\@evenfoot{\@oddfoot}
\makeatother
\fi
\let\maketitlepage\makeagttitle
\let\makeshorttitle\maketitlepage
\let\maketitle\maketitlepage

\newwrite\gtoutfile
\long\gdef\makeheadfile{  %%% start of definition of \makeheadfile
{\def\\{, }\def\s{ }
\immediate\openout\gtoutfile head.xxx
\immediate\write\gtoutfile{To: math@arxiv.org}
\immediate\write\gtoutfile{Subject: put OR rep NNNNN:ppppp}
\immediate\write\gtoutfile{--text follows this line--}
\immediate\write\gtoutfile{Proxy-for: \ifx\theasciiauthors\relax
\theauthors\else\theasciiauthors\fi\s<\ifx\theasciiemail\relax\theemail\else\theasciiemail\fi>}
\immediate\write\gtoutfile{\noexpand\\}
\immediate\write\gtoutfile{Authors: \ifx\theasciiauthors\relax
\theauthors\else\theasciiauthors\fi}
{\def\\{ }\immediate\write\gtoutfile{Title: \ifx\theasciititle\relax
\thetitle\else\theasciititle\fi}}
\immediate\write\gtoutfile{Subj-class: GT or SG, GR etc}
\immediate\write\gtoutfile{MSC-class: \theprimaryclass\ifx\thesecondaryclass\relax\else, \thesecondaryclass\fi}
\immediate\write\gtoutfile{Journal-ref: Algebr. Geom. Topol. \thevolumenumber\s
(\thevolumeyear) \startpage-\finishpage}
\immediate\write\gtoutfile{Comments: Published by Algebraic and
Geometric Topology at}
\immediate\write\gtoutfile{\s\s\s  http://www.maths.warwick.ac.uk/agt/AGTVol\thevolumenumber/agt-\thevolumenumber-\thepapernumber.abs.html}
\immediate\write\gtoutfile{\noexpand\\}
\immediate\write\gtoutfile{}
\ifx\theasciiabstract\relax
\immediate\write\gtoutfile{\theabstract}\else
\immediate\write\gtoutfile{\theasciiabstract}\fi
\immediate\write\gtoutfile{}
\immediate\write\gtoutfile{\noexpand\\}
\immediate\write\gtoutfile{}
\immediate\closeout\gtoutfile}}  %%% end of definition of \makeheadfile

\def\maketitlepage{\makeagttitle\makeheadfile}
\let\makeshorttitle\maketitlepage
\let\maketitle\maketitlepage